\documentclass{amsart}

\usepackage{amssymb}
\usepackage{amsmath}
\usepackage{amsthm}
\usepackage{amsfonts}

\usepackage{a4wide}
\usepackage{graphicx}
\usepackage[usenames,dvipsnames]{xcolor}

\newtheorem{theorem}{Theorem}

\newtheorem{corollary}[theorem]{Corollary}
\newtheorem{lemma}[theorem]{Lemma}

\theoremstyle{definition}

\newtheorem{example}[theorem]{Example}

\newtheorem{counterexample}[theorem]{Counterexample}

\numberwithin{equation}{section}

\begin{document}


\title
{Schur \(\sigma\)-groups with abelian quotient invariants \((9,3)\)}

\author{Daniel C. Mayer}
\address{Naglergasse 53 \\ 8010 Graz \\ Austria}
\email{algebraic.number.theory@algebra.at}
\urladdr{http://www.algebra.at}

\thanks{Research supported by Austrian Science Fund (FWF): P 26008-N25 and Research Executive Agency EUREA}


\date{June 12, 2020}


\begin{abstract}
By the construction of suitable non-metabelian Schur \(\sigma\)-groups \(S\) of type \((9,3)\)
with order \(\#(S)=3^{21}\) and nilpotency class \(\mathrm{cl}(S)=9\),
evidence is provided of a new class of imaginary quadratic fields \(K\)
with \(3\)-class group \(\mathrm{Cl}_3(K)\simeq C_9\times C_3\)
and punctured principalization type \(\varkappa\sim (1,4,4;4)\)
whose \(3\)-class field tower consists of precisely three stages.
In contrast, previous examples of three-stage towers were associated with
\(\varkappa\in\lbrace (1,1,2;2), (1,1,2;3), (1,1,4;2), (1,2,3;1)\rbrace\),
\(\#(S)=3^9\) and \(\mathrm{cl}(S)=5\).
\end{abstract}

\maketitle


\section{Introduction}
\label{s:Intro}

\noindent
The aim of this article is
the construction of Schur \(\sigma\)-groups \(S\)
with non-elementary bicyclic derived quotient \(S/S^\prime\)
isomorphic to the product \(C_9\times C_3\) of cyclic groups
and the realization of such groups as automorphism groups
\(\mathrm{G}_3^\infty(K):=\mathrm{Gal}(\mathrm{F}_3^\infty(K)/K)\)
of maximal unramified pro-\(3\)-extensions \(\mathrm{F}_3^\infty(K)\)
of imaginary quadratic fields \(K=\mathbb{Q}(\sqrt{d})\), \(d<0\),
with \(3\)-class group \(\mathrm{Cl}_3(K)\) of type \((9,3)\).
Based on the fundamental Shafarevich theorem
\cite{Sh,Ma2015d}
concerning the relation rank
of the pro-\(p\)-group \(\mathfrak{G}:=\mathrm{G}_p^\infty(K)\)
of the \(p\)-class field tower
\(K\le\mathrm{F}_p^1(K)\le\mathrm{F}_p^2(K)\le\ldots\le\mathrm{F}_p^\infty(K)\)
of an algebraic number field \(K\),
Koch and Venkov
\cite{KoVe}
have shown that the Galois group \(\mathfrak{G}\)
must be a \textit{Schur \(\sigma\)-group}
\cite{Ag}
with balanced presentation and generator inverting \(\sigma\)-automorphism,
when \(K\) is an imaginary quadratic field
and \(p\) is an odd prime number.
In this paper, the focus is on the smallest odd prime \(p=3\).

In comparison to \(3\)-groups \(S\)
with elementary abelianization \(S/S^\prime\simeq C_3\times C_3\) of rank two,
only modest progress has been made in the classification of
\(3\)-groups \(S\) with non-elementary \(S/S^\prime\simeq C_9\times C_3\).
With respect to the complexity of the structure,
it turns out that \textit{for both types} of derived quotients \(S/S^\prime\),
there appear the same \textit{three phenomena} with entirely different characteristics:
\begin{itemize}
\item
\textbf{Phenomenon 1.}
Existence of \textit{finitely} many \textit{metabelian} Schur \(\sigma\)-groups \(S\)
with nilpotency class \(\mathrm{cl}(S)=3\),
corresponding to \(3\)-class field towers \(\mathrm{F}_3^2(K)=\mathrm{F}_3^\infty(K)\) with two stages,
enjoying a dense population \(S\simeq\mathrm{G}_3^\infty(K)\) by imaginary quadratic fields \(K\).
\item
\textbf{Phenomenon 2.}
Existence of an \textit{infinitude} of \textit{non-metabelian} Schur \(\sigma\)-groups \(S\)
with unbounded nilpotency class \(\mathrm{cl}(S)\ge 5\),
pairwise non-isomorphic second derived quotients \(S/S^{\prime\prime}\),
and \textit{fixed derived length} \(\mathrm{dl}(S)=3\),
corresponding to \(3\)-class field towers \(\mathrm{F}_3^3(K)=\mathrm{F}_3^\infty(K)\) with precisely three stages,
populated moderately but significantly by imaginary quadratic fields \(K\).
\item
\textbf{Phenomenon 3.}
Existence of \textit{infinitely} many \textit{multiplets} \((S_i)_{1\le i\le m}\) with multiplicity \(m\)
of \textit{non-metabelian} Schur \(\sigma\)-groups,
sharing a common metabelianization \(S_i/S_i^{\prime\prime}\),
with unbounded nilpotency class \(\mathrm{cl}(S_i)\ge 9\)
and \textit{unbounded derived length} \(\mathrm{dl}(S_i)\ge 3\),
corresponding to \(3\)-class field towers \(\mathrm{F}_3^\infty(K)\) with at least three stages,
populated only sparsely by imaginary quadratic fields \(K\).
\end{itemize}

\noindent
Since basic features of metabelian Schur \(\sigma\)-groups \(S\),
without explicit mention of the relation rank,
have been investigated completely by Scholz and Taussky
\cite{SoTa}[Section D]
--- in the notation of the SmallGroups database
\cite{BEO},
there are only two such groups \(S\in\lbrace\langle 243,5\rangle,\langle 243,7\rangle\rbrace\)
with \(S/S^\prime\simeq C_3\times C_3\),
and also two \(S\in\lbrace\langle 729,14\rangle,\langle 729,15\rangle\rbrace\)
with \(S/S^\prime\simeq C_9\times C_3\) ---
and since Phenomenon 2 with bijective correspondence
between non-metabelian Schur \(\sigma\)-groups \(S\) and second derived quotients \(S/S^{\prime\prime}\) 
has been treated exhaustively in a series of recent publications
\cite{BuMa,Ma2015a,Ma2015c,Ma2018b},
the focus in this paper is on Phenomenon 3
for non-elementary abelianization.

It must be pointed out that useful foundations,
though restricted to metabelianizations and without consideration of relation ranks,
for Phenomenon 2 have been established with remarkable presentiment by Scholz and Taussky in
\cite{SoTa}[Section E]
and for Phenomenon 3 in
\cite{SoTa}[Section F].
These foundations were confirmed and extended with modern techniques by Nebelung in
\cite{Ne1,Ne2}.


\section{Abelian quotient invariants of second order}
\label{s:AQI}

\noindent
In order to identify
non-metabelian Schur \(\sigma\)-groups \(S\)
with nilpotency class \(\mathrm{cl}(S)=9\)
and derived length \(\mathrm{dl}(S)=3\),
their \textit{logarithmic abelian quotient invariants} (AQI) \textit{of second order}
\begin{equation}
\label{eqn:AQI}
\tau^{(2)}(S)=\lbrack 21;(\tau_0;2^21^3,T_1);(2^21;2^21^3,T_2),(2^21;2^21^3,T_3),(2^21;2^21^3,T_4)\rbrack,
\end{equation}
i.e., the AQI of subgroups \(U\) of index \((S:U)=9\), can be used,
where \(\tau_0\) is either \textit{homo}cyclic of type \((2^3)\) or \textit{hetero}cyclic of type \((321)\)
and the possible structures of thirteen components of \(T_j\) with \(1\le j\le 4\)
are abbreviated by the \textit{symbols} in Table
\ref{tbl:Symbols}.


\renewcommand{\arraystretch}{1.2}

\begin{table}[ht]
\caption{Symbols for components of AQI of second order}
\label{tbl:Symbols}
\begin{center}
\begin{tabular}{|c||c|c|}
\hline
 Symbol          & \multicolumn{2}{|c|}{Structure} \\
\hline
 \(\alpha\)      & \((321)^3\)    & \((32)^9\)   \\
 \(\beta\)       & \((321)^3\)    & \((2^21)^9\) \\
 \(\gamma\)      & \((31^3)^3\)   & \((32)^9\)   \\
 \(\delta\)      & \((2^21^2)^3\) & \((2^21)^9\) \\
 \(\varepsilon\) & \((321^2)^3\)  & \((2^21)^9\) \\
 \(\zeta\)       & \((321^3)^3\)  & \((2^21)^9\) \\
 \(\eta\)        & \((321^2)^3\)  & \((32)^9\)   \\
 \(\vartheta\)   & \((2^31)^3\)   & \((2^21)^9\) \\
 \(\xi\)         & \((41^4)^3\)   & \((32)^9\)   \\
\hline
\end{tabular}
\end{center}
\end{table}


\noindent
The \(18\) relevant ancestors of the desired Schur \(\sigma\)-groups \(S\)
are immediate descendants \(P_8-\#4;\ell\) with \(1\le\ell\le 72\) and step size \(s=4\)
of the fork vertex \(P_8:=\langle 6561,165\rangle=\langle 729,10\rangle-\#2;2\).
Their second AQI are listed together with information on variants in Table
\ref{tbl:HigherAQI}.

\newpage


\renewcommand{\arraystretch}{1.2}

\begin{table}[ht]
\caption{AQI of second order represented as punctured quartets}
\label{tbl:HigherAQI}
\begin{center}
\begin{tabular}{|r|cc|c|c|c|c|}
\hline
 Identifier & \multicolumn{5}{|c|}{\(\tau^{(2)}(S)\)} & \\
 \(\ell\)   & \(\tau_0\) & \(T_1\)         & \(T_2\)    & \(T_3\)    & \(T_4\)    & Variant \\
\hline
  \(2\)     & \((2^3)\)  & \(\zeta\)       & \(\alpha\) & \(\alpha\) & \(\alpha\) & Triplet \\
  \(4\)     & \((2^3)\)  & \(\varepsilon\) & \(\alpha\) & \(\alpha\) & \(\delta\) & Triplet \\
  \(6\)     & \((2^3)\)  & \(\varepsilon\) & \(\alpha\) & \(\alpha\) & \(\delta\) & Singlet \\
  \(8\)     & \((2^3)\)  & \(\zeta\)       & \(\alpha\) & \(\beta\)  & \(\beta\)  & \\
  \(9\)     & \((2^3)\)  & \(\zeta\)       & \(\gamma\) & \(\gamma\) & \(\gamma\) & \\
 \(11\)     & \((2^3)\)  & \(\varepsilon\) & \(\alpha\) & \(\alpha\) & \(\gamma\) & Triplet \\
 \(12\)     & \((2^3)\)  & \(\varepsilon\) & \(\alpha\) & \(\beta\)  & \(\gamma\) & Triplet \\
 \(14\)     & \((2^3)\)  & \(\varepsilon\) & \(\alpha\) & \(\alpha\) & \(\gamma\) & Singlet \\
 \(15\)     & \((2^3)\)  & \(\varepsilon\) & \(\alpha\) & \(\beta\)  & \(\gamma\) & Singlet \\
\hline
 \(17\)     & \((321)\)  & \(\eta\)        & \(\alpha\) & \(\beta\)  & \(\gamma\) & Singlet \\
 \(19\)     & \((321)\)  & \(\eta\)        & \(\alpha\) & \(\alpha\) & \(\gamma\) & Singlet \\
 \(21\)     & \((321)\)  & \(\xi\)         & \(\alpha\) & \(\alpha\) & \(\beta\)  & Triplet \\
 \(23\)     & \((321)\)  & \(\eta\)        & \(\alpha\) & \(\alpha\) & \(\gamma\) & Triplet \\
 \(24\)     & \((321)\)  & \(\eta\)        & \(\alpha\) & \(\beta\)  & \(\gamma\) & Singlet \\
 \(26\)     & \((321)\)  & \(\eta\)        & \(\alpha\) & \(\beta\)  & \(\delta\) & \\
 \(27\)     & \((321)\)  & \(\eta\)        & \(\alpha\) & \(\alpha\) & \(\gamma\) & Triplet \\
 \(29\)     & \((321)\)  & \(\xi\)         & \(\alpha\) & \(\alpha\) & \(\alpha\) & Singlet \\
 \(30\)     & \((321)\)  & \(\xi\)         & \(\gamma\) & \(\gamma\) & \(\delta\) & \\
\hline
\end{tabular}
\end{center}
\end{table}


\section{Main Theorem for non-elementary \(S/S^\prime\simeq C_9\times C_3\)}
\label{s:Main93}

\noindent
The following theorem provides evidence of
a new class of algebraic number fields of type \((9,3)\)
whose \(3\)-class field tower consists of exactly three stages.

\begin{theorem}
\label{thm:ThreeStageTower93}
An imaginary quadratic field \(K=\mathbb{Q}(\sqrt{d})\)
with non-elementary \(3\)-class group \(\mathrm{Cl}_3(K)\simeq C_9\times C_3\) of rank two,
punctured capitulation type \(\mathrm{B}.18\), i.e.
\[\varkappa(K)\sim (144;4),\]
and abelian type invariants \(\tau^{(2)}(K)\) of second order of the shape in Formula
\eqref{eqn:AQI}
with either
\begin{equation}
\label{eqn:AQI1}
\tau_0=2^3, \quad T_1=\zeta, \quad T_2=\alpha, \quad T_3=\alpha, \quad T_4=\alpha
\end{equation}
or
\begin{equation}
\label{eqn:AQI2}
\tau_0=2^3, \quad T_1=\varepsilon, \quad T_2=\alpha, \quad T_3=\alpha, \quad T_4=\delta
\end{equation}
or
\begin{equation}
\label{eqn:AQI3}
\tau_0=2^3, \quad T_1=\varepsilon, \quad T_2=\alpha, \quad T_3=\alpha, \quad T_4=\gamma
\end{equation}
or
\begin{equation}
\label{eqn:AQI4}
\tau_0=2^3, \quad T_1=\varepsilon, \quad T_2=\alpha, \quad T_3=\beta, \quad T_4=\gamma
\end{equation}
or
\begin{equation}
\label{eqn:AQI5}
\tau_0=321, \quad T_1=\eta, \quad T_2=\alpha, \quad T_3=\alpha, \quad T_4=\gamma
\end{equation}
or
\begin{equation}
\label{eqn:AQI6}
\tau_0=321, \quad T_1=\xi, \quad T_2=\alpha, \quad T_3=\alpha, \quad T_4=\alpha
\end{equation}
possesses a finite \(3\)-class field tower
\[K=\mathrm{F}_3^0(K)<\mathrm{F}_3^1(K)<\mathrm{F}_3^2(K)<\mathrm{F}_3^3(K)=\mathrm{F}_3^\infty(K)\]
with precise length \(\ell_3(K)=3\).
\end{theorem}


In the following corollary,
Theorem
\ref{thm:ThreeStageTower93}
is supplemented by information on the Galois group
\(\mathfrak{G}:=\mathrm{Gal}(\mathrm{F}_3^\infty(K)/K)\)
and its metabelianization
\(\mathfrak{G}/\mathfrak{G}^{\prime\prime}\simeq\mathrm{Gal}(\mathrm{F}_3^2(K)/K)\).

\begin{corollary}
\label{cor:ThreeStageTower93}
Let \(K\) be a field with properties as in the assumptions of Theorem
\ref{thm:ThreeStageTower93}.
Then the automorphism group \(\mathrm{Gal}(\mathrm{F}_3^\infty(K)/K)\)
of the full \(3\)-class field tower of \(K\) is a non-metabelian Schur \(\sigma\)-group
\cite{Ag,KoVe}
with derived length \(3\), order \(3^{21}\) and nilpotency class \(9\).
The second \(3\)-class group \(\mathrm{Gal}(\mathrm{F}_3^2(K)/K)\) of \(K\)
\cite{Ma2012}
is a metabelian \(\sigma\)-group of order \(3^{10}\) and nilpotency class \(5\).
\end{corollary}


\begin{example}
\label{exm:ThreeStageTower93}
The quadratic fields \(K\) with fundamental discriminants \\
\(d_K\in\lbrace -1\,103\,784, -1\,356\,215\rbrace\) and \(\tau^{(2)}(K)\) in
\eqref{eqn:AQI1}, \\
respectively \(d_K=-518\,835\) and \(\tau^{(2)}(K)\) in
\eqref{eqn:AQI2}, \\
respectively \(d_K=-761\,855\) and \(\tau^{(2)}(K)\) in
\eqref{eqn:AQI3}, \\
respectively
\(d_K\in\lbrace -553\,807, -763\,972, -876\,948, -1\,100\,315, -1\,407\,379, -1\,677\,327, -1\,727\,187\rbrace\)
and \(\tau^{(2)}(K)\) in
\eqref{eqn:AQI4}, \\
respectively \(d_K\in\lbrace -194\,703, -494\,771\rbrace\) and \(\tau^{(2)}(K)\) in
\eqref{eqn:AQI5}, \\
respectively \(d_K=-150\,319\) and \(\tau^{(2)}(K)\) in
\eqref{eqn:AQI6}, \\
satisfy the assumptions of Theorem
\ref{thm:ThreeStageTower93}
with punctured capitulation type \(\varkappa(K)\sim (144;4)\).
Consequently,
each of them is an example of a field possessing a
\(3\)-class field tower with exactly three stages of relative degrees
\(\lbrack\mathrm{F}_3^3(K):\mathrm{F}_3^2(K)\rbrack=3^{11}\),
\(\lbrack\mathrm{F}_3^2(K):\mathrm{F}_3^1(K)\rbrack=3^7\),
\(\lbrack\mathrm{F}_3^1(K):\mathrm{F}_3^0(K)\rbrack=3^3\),
and Galois group \(\mathrm{Gal}(\mathrm{F}_3^\infty(K)/K)\) of order \(3^{21}\).
\end{example}


\begin{counterexample}
\label{cex:ThreeStageTower93}
Unfortunately, the quadratic fields \(K\) with fundamental discriminants
\(d_K\in\lbrace -294\,983, -389\,435\rbrace\)
and punctured capitulation type \(\varkappa(K)\sim (144;4)\)
do not satisfy the assumptions of Theorem
\ref{thm:ThreeStageTower93},
since \(\tau^{(2)}(K)\)
is of the shape in Formula
\eqref{eqn:AQI}
with
\[\tau_0=2^3, \quad T_1=\zeta, \quad T_2=\varepsilon, \quad T_3=\eta, \quad T_4=\vartheta.\]
In these cases, the number of stages is unknown
and may even be infinite.
\end{counterexample}



\section{Proof of the Main Theorem for \(S/S^\prime\simeq C_9\times C_3\)}
\label{s:Proof93}

\noindent
Let \(P_8:=\langle 6561,165\rangle=\langle 729,10\rangle-\#2;2\) (identifier in
\cite{BEO})
denote the common fork
of the root paths of all finite \(3\)-groups \(G\)
with non-elementary abelianization \(G/G^\prime\simeq C_9\times C_3\) of rank two,
punctured transfer kernel type \(\mathrm{B}.18\),
\[\varkappa(G)\sim (144;4),\]
and abelian quotient invariants of first order
\[\tau^{(1)}(G)=\Bigl((9,3);\tau_0,(9,3,3)^3\Bigr)\]
with \(\tau_0\in\lbrace (9,9,9), (27,9,3)\rbrace\).


\begin{lemma}
\label{lem:Type1}
There exist precisely \(81\) Schur \(\sigma\)-groups \(S\)
with non-elementary abelianization \(S/S^\prime\simeq C_9\times C_3\) of rank two,
punctured transfer kernel type \(\mathrm{B}.18\),
\(\varkappa(S)\sim (144;4)\),
and abelian quotient invariants \(\tau^{(2)}(S)\) of second order in Formula
\eqref{eqn:AQI1}.
Their relative identifiers
\cite{MAGMA}
with respect to \(P_8\) are uniformly given in the shape
\[P_8-\#4;\ell-\#2;k-\#4;j-\#1;i-\#2;h\]
with a single possible pair \((\ell,k)\), \(1\le\ell\le 72\), \(1\le k\le 41\),
namely \((\ell,k)=(2,5)\),
and \(27\) possible pairs \((j,i)\), \(1\le j\le 27\), \(1\le i\le 5\).
Here, \(k\) is determined uniquely as a function \(k=k(\ell)\) of \(\ell\),
\(j\) runs through all possible values, 
\(i\) is determined uniquely as a function \(i=i(j)\) of \(j\),
and \(1\le h\le 3\).
The metabelianization \(M:=S/S^{\prime\prime}\) is \(M\simeq P_8-\#2;83\).
The groups \(S\) are of order \(3^{21}\), class \(9\), coclass \(12\),
and their common metabelianization \(M\) is of order \(3^{10}\), class \(5\), coclass \(5\).
\end{lemma}


\begin{lemma}
\label{lem:Type2}
There exist precisely \(108=81+27\) Schur \(\sigma\)-groups \(S\)
with non-elementary abelianization \(S/S^\prime\simeq C_9\times C_3\) of rank two,
punctured transfer kernel type \(\mathrm{B}.18\),
\(\varkappa(S)\sim (144;4)\),
and abelian quotient invariants \(\tau^{(2)}(S)\) of second order in Formula
\eqref{eqn:AQI2}.
Their relative identifiers
\cite{MAGMA}
with respect to \(P_8\) are uniformly given in the shape
\[P_8-\#4;\ell-\#2;k-\#4;j-\#1;i-\#2;h\]
with two possible pairs \((\ell,k)\), \(1\le\ell\le 72\), \(1\le k\le 41\),
namely \((\ell,k)\in\lbrace (4,41), (6,41)\rbrace\),
and, for each of them, \(27\) possible pairs \((j,i)\), \(1\le j\le 27\), \(1\le i\le 5\).
Here, \(k\) is determined uniquely as a function \(k=k(\ell)\) of \(\ell\),
\(j\) runs through all possible values, 
\(i\) is determined uniquely as a function \(i=i(j)\) of \(j\),
and \(1\le h\le 3\) for the first pair, \(h=1\) for the second pair.
The metabelianization \(M:=S/S^{\prime\prime}\) is \(M\simeq P_8-\#2;83\).
The groups \(S\) are of order \(3^{21}\), class \(9\), coclass \(12\),
and their common metabelianization \(M\) is of order \(3^{10}\), class \(5\), coclass \(5\).
\end{lemma}


\begin{lemma}
\label{lem:Type3}
There exist precisely \(108=81+27\) Schur \(\sigma\)-groups \(S\)
with non-elementary abelianization \(S/S^\prime\simeq C_9\times C_3\) of rank two,
punctured transfer kernel type \(\mathrm{B}.18\),
\(\varkappa(S)\sim (144;4)\),
and abelian quotient invariants \(\tau^{(2)}(S)\) of second order in Formula
\eqref{eqn:AQI3}.
Their relative identifiers
\cite{MAGMA}
with respect to \(P_8\) are uniformly given in the shape
\[P_8-\#4;\ell-\#2;k-\#4;j-\#1;i-\#2;h\]
with two possible pairs \((\ell,k)\), \(1\le\ell\le 72\), \(1\le k\le 41\),
namely \((\ell,k)\in\lbrace (11,41), (14,11)\rbrace\),
and, for each of them, \(27\) possible pairs \((j,i)\), \(1\le j\le 27\), \(1\le i\le 5\).
Here, \(k\) is determined uniquely as a function \(k=k(\ell)\) of \(\ell\),
\(j\) runs through all possible values, 
\(i\) is determined uniquely as a function \(i=i(j)\) of \(j\),
and \(1\le h\le 3\) for the first pair, \(h=1\) for the second pair.
The metabelianization \(M:=S/S^{\prime\prime}\) is \(M\simeq P_8-\#2;83\).
The groups \(S\) are of order \(3^{21}\), class \(9\), coclass \(12\),
and their common metabelianization \(M\) is of order \(3^{10}\), class \(5\), coclass \(5\).
\end{lemma}


\begin{lemma}
\label{lem:Type4}
There exist precisely \(108=81+27\) Schur \(\sigma\)-groups \(S\)
with non-elementary abelianization \(S/S^\prime\simeq C_9\times C_3\) of rank two,
punctured transfer kernel type \(\mathrm{B}.18\),
\(\varkappa(S)\sim (144;4)\),
and abelian quotient invariants \(\tau^{(2)}(S)\) of second order in Formula
\eqref{eqn:AQI4}.
Their relative identifiers
\cite{MAGMA}
with respect to \(P_8\) are uniformly given in the shape
\[P_8-\#4;\ell-\#2;k-\#4;j-\#1;i-\#2;h\]
with two possible pairs \((\ell,k)\), \(1\le\ell\le 72\), \(1\le k\le 41\),
namely \((\ell,k)\in\lbrace (12,41), (15,11)\rbrace\),
and, for each of them, \(27\) possible pairs \((j,i)\), \(1\le j\le 27\), \(1\le i\le 5\).
Here, \(k\) is determined uniquely as a function \(k=k(\ell)\) of \(\ell\),
\(j\) runs through all possible values, 
\(i\) is determined uniquely as a function \(i=i(j)\) of \(j\),
and \(1\le h\le 3\) for the first pair, \(h=1\) for the second pair.
The metabelianization \(M:=S/S^{\prime\prime}\) is \(M\simeq P_8-\#2;83\).
The groups \(S\) are of order \(3^{21}\), class \(9\), coclass \(12\),
and their common metabelianization \(M\) is of order \(3^{10}\), class \(5\), coclass \(5\).
\end{lemma}


\begin{lemma}
\label{lem:Type5}
There exist precisely \(189=27+81+81\) Schur \(\sigma\)-groups \(S\)
with non-elementary abelianization \(S/S^\prime\simeq C_9\times C_3\) of rank two,
punctured transfer kernel type \(\mathrm{B}.18\),
\(\varkappa(S)\sim (144;4)\),
and abelian quotient invariants \(\tau^{(2)}(S)\) of second order in Formula
\eqref{eqn:AQI5}.
Their relative identifiers
\cite{MAGMA}
with respect to \(P_8\) are uniformly given in the shape
\[P_8-\#4;\ell-\#2;k-\#4;j-\#1;i-\#2;h\]
with three possible pairs \((\ell,k)\), \(1\le\ell\le 72\), \(1\le k\le 41\),
namely \((\ell,k)\in\lbrace (19,41), (23,14), (27,11)\rbrace\),
and, for each of them, \(27\) possible pairs \((j,i)\), \(1\le j\le 27\), \(1\le i\le 5\).
Here, \(k\) is determined uniquely as a function \(k=k(\ell)\) of \(\ell\),
\(j\) runs through all possible values, 
\(i\) is determined uniquely as a function \(i=i(j)\) of \(j\),
and \(h=1\) for the first pair, \(1\le h\le 3\) for the second and third pair.
The metabelianization \(M:=S/S^{\prime\prime}\) is \(M\simeq P_8-\#2;85\).
The groups \(S\) are of order \(3^{21}\), class \(9\), coclass \(12\),
and their common metabelianization \(M\) is of order \(3^{10}\), class \(5\), coclass \(5\).
\end{lemma}


\begin{lemma}
\label{lem:Type6}
There exist precisely \(27\) Schur \(\sigma\)-groups \(S\)
with non-elementary abelianization \(S/S^\prime\simeq C_9\times C_3\) of rank two,
punctured transfer kernel type \(\mathrm{B}.18\),
\(\varkappa(S)\sim (144;4)\),
and abelian quotient invariants \(\tau^{(2)}(S)\) of second order in Formula
\eqref{eqn:AQI6}.
Their relative identifiers
\cite{MAGMA}
with respect to \(P_8\) are uniformly given in the shape
\[P_8-\#4;\ell-\#2;k-\#4;j-\#1;i-\#2;h\]
with a single possible pair \((\ell,k)\), \(1\le\ell\le 72\), \(1\le k\le 41\),
namely \((\ell,k)=(29,41)\),
and \(27\) possible pairs \((j,i)\), \(1\le j\le 27\), \(1\le i\le 5\).
Here, \(k\) is determined uniquely as a function \(k=k(\ell)\) of \(\ell\),
\(j\) runs through all possible values, 
\(i\) is determined uniquely as a function \(i=i(j)\) of \(j\),
and \(h=1\).
The metabelianization \(M:=S/S^{\prime\prime}\) is \(M\simeq P_8-\#2;85\).
The groups \(S\) are of order \(3^{21}\), class \(9\), coclass \(12\),
and their common metabelianization \(M\) is of order \(3^{10}\), class \(5\), coclass \(5\).
\end{lemma}


The lemmas have been proven with the aid of the
computational algebra system MAGMA
\cite{MAGMA},
developed at the University of Sydney
under supervision of John Cannon
\cite{BCP,BCFS}.

For the proof of Theorem
\ref{thm:ThreeStageTower93}
and Corollary
\ref{cor:ThreeStageTower93},
the invariants \(\varkappa(K)\), \(\tau^{(1)}(K)\) and \(\tau^{(2)}(K)\)
of \(3\)-class groups associated with algebraic number fields \(K\)
are first translated to invariants \(\varkappa(G)\), \(\tau^{(1)}(G)\) and \(\tau^{(2)}(G)\)
of finite \(3\)-groups, according to the Artin reciprocity law of class field theory
\cite{Ar1,Ar2}.

Then the group theoretic invariants are used as a search pattern
for Schur \(\sigma\)-groups, which are required as Galois group
\(\mathrm{Gal}(\mathrm{F}_3^\infty(K)/K)\)
of imaginary quadratic fields \(K\)
\cite{Ag,KoVe}.

Since Schur \(\sigma\)-groups have extremal root paths,
the search can be restricted to \(\sigma\)-groups connected by edges with maximal step size
equal to the nuclear rank of the respective parent.

The nuclear rank of the fork \(P_8\) is \(4\),
and the various search patterns \(\varkappa(G)\), \(\tau^{(1)}(G)\) and \(\tau^{(2)}(G)\)
are hit by the respective descendants \(P_8-\#4;\ell\)
in the Lemmas
\ref{lem:Type1}
\(\ldots\)
\ref{lem:Type6}.

The pruned trees of \(\sigma\)-descendants are finite
and terminate in the Schur \(\sigma\)-groups
\(P_8-\#4;\ell-\#2;k-\#4;j-\#1;i-\#2;h\)
of order \(3^{21}\) as leaves. QED.


\section{Main Theorem for elementary \(S/S^\prime\simeq C_3\times C_3\)}
\label{s:Main33}

\noindent
In order to enable comparison with elementary abelianization,
the following theorem provides evidence of
a new class of algebraic number fields of type \((3,3)\)
whose \(3\)-class field tower consists of exactly three stages.

\begin{theorem}
\label{thm:ThreeStageTower33}
An imaginary quadratic field \(K=\mathbb{Q}(\sqrt{d})\)
with elementary \(3\)-class group \(\mathrm{Cl}_3(K)\) of rank two,
capitulation type either \(\mathrm{F}.12\) or \(\mathrm{F}.13\)
\cite{SoTa,Ma1991,Ma2010},
i.e.
\[\text{either} \quad\varkappa(K)\sim (2114)\quad \text{or} \quad\varkappa(K)\sim (2141),\]
and abelian type invariants of second order either
\[\tau^{(2)}(K)=\Bigl((3,3);\lbrack (27,9);(9,9,9,3),(27,3,3,3)^3\rbrack^2,\lbrack (3,3,3);(9,9,9,3),(9,9,3)^3,(9,3,3)^9\rbrack^2\Bigr)\]
possesses a finite \(3\)-class field tower
\[K=\mathrm{F}_3^0(K)<\mathrm{F}_3^1(K)<\mathrm{F}_3^2(K)<\mathrm{F}_3^3(K)=\mathrm{F}_3^\infty(K)\]
with precise length \(\ell_3(K)=3\).
\end{theorem}

Both capitulation types, \(\varkappa(K)\sim (2114)\) and \(\varkappa(K)\sim (2141)\),
contain a transposition, \(1\mapsto 2\) and \(2\mapsto 1\),
but the former has a fixed point \(4\mapsto 4\),
whereas the latter has no fixed points.


In the following corollary,
Theorem
\ref{thm:ThreeStageTower33}
is supplemented by information on the Galois group
\(\mathfrak{G}:=\mathrm{Gal}(\mathrm{F}_3^\infty(K)/K)\)
and its metabelianization
\(\mathfrak{G}/\mathfrak{G}^{\prime\prime}\simeq\mathrm{Gal}(\mathrm{F}_3^2(K)/K)\).

\begin{corollary}
\label{cor:ThreeStageTower33}
Let \(K\) be a field with properties as in the assumptions of Theorem
\ref{thm:ThreeStageTower33}.
Then the automorphism group \(\mathrm{Gal}(\mathrm{F}_3^\infty(K)/K)\)
of the full \(3\)-class field tower of \(K\) is a non-metabelian Schur \(\sigma\)-group
\cite{Ag,KoVe}
with derived length \(3\), order \(3^{20}\) and coclass \(11\).
The second \(3\)-class group \(\mathrm{Gal}(\mathrm{F}_3^2(K)/K)\) of \(K\)
\cite{Ma2012}
is a metabelian \(\sigma\)-group of order \(3^9\) and coclass \(4\).
\end{corollary}


\begin{example}
\label{exm:ThreeStageTower33}
The quadratic fields \(K\) with fundamental discriminants
\(d_K\in\lbrace -291\,220, -633\,352\), \(-874\,680\rbrace\),
respectively \(d_K=-731\,867\),
satisfy the assumptions of Theorem
\ref{thm:ThreeStageTower33}
with capitulation type \(\varkappa(K)\sim (2114)\),
respectively \(\varkappa(K)\sim (2141)\).
Consequently,
each of them is an example of a field possessing a
\(3\)-class field tower with exactly three stages of relative degrees
\(\lbrack\mathrm{F}_3^3(K):\mathrm{F}_3^2(K)\rbrack=3^{11}\),
\(\lbrack\mathrm{F}_3^2(K):\mathrm{F}_3^1(K)\rbrack=3^7\),
\(\lbrack\mathrm{F}_3^1(K):\mathrm{F}_3^0(K)\rbrack=3^2\),
and Galois group \(\mathrm{Gal}(\mathrm{F}_3^\infty(K)/K)\) of order \(3^{20}\).
\end{example}


\begin{counterexample}
\label{cex:ThreeStageTower33}
Unfortunately, the quadratic fields \(K\) with fundamental discriminants
\(d_K\in\lbrace -507\,140, -578\,847\rbrace\)
and capitulation type \(\varkappa(K)\sim (2114)\),
respectively
\(d_K\in\lbrace -167\,064\),
\( -296\,407, -317\,747, -401\,603, -588\,027, -591\,412, -803\,591, -835\,707, -999\,704\rbrace\),
which share the common capitulation type \(\varkappa(K)\sim (2141)\),
almost but not completely satisfy the assumptions of Theorem
\ref{thm:ThreeStageTower33},
since \(\tau^{(2)}(K)\)
contains a single inadequate component
\[\lbrack (27,9);(9,9,9,3),(81,3,3,3)^3\rbrack.\]
In these cases, the number of stages can only be estimated by
\(3\le\ell_3(K)\le 4\)
and \(\mathrm{Gal}(\mathrm{F}_3^\infty(K)/K)\)
is of order either \(3^{20}\) with derived length \(3\)
or \(3^{23}\) with derived length \(4\).
\end{counterexample}



\section{Proof of the Main Theorem for \(S/S^\prime\simeq C_3\times C_3\)}
\label{s:Proof33}

\noindent
Let \(P_7:=\langle 2187,64\rangle=\langle 243,3\rangle-\#2;1\) (identifier in
\cite{BEO})
denote the common fork
of the root paths of all finite \(3\)-groups \(G\)
with elementary abelianization \(G/G^\prime\) of rank two,
transfer kernel type
\[\text{either} \quad\varkappa(G)\sim (2114)\quad \text{or} \quad\varkappa(G)\sim (2141),\]
and abelian quotient invariants of first order
\[\tau^{(1)}(G)=\Bigl((3,3);(27,9)^2,(3,3,3)^2\Bigr).\]


\begin{lemma}
\label{lem:TypeF12}
There exist precisely \(216=8\cdot 27\) Schur \(\sigma\)-groups \(S\)
with elementary abelianization \(S/S^\prime\) of rank two,
transfer kernel type \(\mathrm{F}.12\),
\(\varkappa(S)\sim (2114)\),
and abelian quotient invariants of second order
\[\tau^{(2)}(S)=\Bigl((3,3);\lbrack (27,9);(9,9,9,3),(27,3,3,3)^3\rbrack^2,\lbrack (3,3,3);(9,9,9,3),(9,9,3)^3,(9,3,3)^9\rbrack^2\Bigr).\]
Their relative identifiers
\cite{MAGMA}
with respect to \(P\) are uniformly given in the shape
\[P_7-\#4;\ell-\#2;k-\#4;j-\#1;i-\#2;1\]
with eight possible pairs \((\ell,k)\), \(1\le\ell\le 198\), \(1\le k\le 41\),
and, for each of them, \(27\) possible pairs \((j,i)\), \(1\le j\le 27\), \(1\le i\le 5\).
Here, \(k\) is determined uniquely as a function \(k=k(\ell)\) of \(\ell\),
\(j\) runs through all possible values, and
\(i\) is determined uniquely as a function \(i=i(j)\) of \(j\).
According to the metabelianization \(M:=S/S^{\prime\prime}\), there are four cases,
\begin{enumerate}
\item
\(M\simeq P_7-\#2;43\) and \((\ell,k)\in\lbrace (125,11),(143,14)\rbrace\),
\item
\(M\simeq P_7-\#2;46\) and \((\ell,k)\in\lbrace (130,1),(146,13)\rbrace\),
\item
\(M\simeq P_7-\#2;51\) and \((\ell,k)\in\lbrace (113,32),(126,41)\rbrace\),
\item
\(M\simeq P_7-\#2;53\) and \((\ell,k)\in\lbrace (116,31),(118,41)\rbrace\).
\end{enumerate}
The groups \(S\) are of order \(3^{20}\), class \(9\), coclass \(11\),
and their metabelianizations \(M\) are of order \(3^9\), class \(5\), coclass \(4\).
\end{lemma}


\begin{lemma}
\label{lem:TypeF13}
There exist precisely \(216=8\cdot 27\) Schur \(\sigma\)-groups \(S\)
with elementary abelianization \(S/S^\prime\) of rank two,
transfer kernel type \(\mathrm{F}.13\),
\(\varkappa(S)\sim (2141)\),
and abelian quotient invariants of second order
\[\tau^{(2)}(S)=\Bigl((3,3);\lbrack (27,9);(9,9,9,3),(27,3,3,3)^3\rbrack^2,\lbrack (3,3,3);(9,9,9,3),(9,9,3)^3,(9,3,3)^9\rbrack^2\Bigr).\]
Their relative identifiers
\cite{MAGMA}
with respect to \(P\) are uniformly given in the shape
\[P_7-\#4;\ell-\#2;k-\#4;j-\#1;i-\#2;1\]
with eight possible pairs \((\ell,k)\), \(1\le\ell\le 198\), \(1\le k\le 41\),
and, for each of them, \(27\) possible pairs \((j,i)\), \(1\le j\le 27\), \(1\le i\le 5\).
Here, \(k\) is determined uniquely as a function \(k=k(\ell)\) of \(\ell\),
\(j\) runs through all possible values, and
\(i\) is determined uniquely as a function \(i=i(j)\) of \(j\).
According to the metabelianization \(M:=S/S^{\prime\prime}\), there are four cases,
\begin{enumerate}
\item
\(M\simeq P_7-\#2;41\) and \((\ell,k)\in\lbrace (177,1),(185,41)\rbrace\),
\item
\(M\simeq P_7-\#2;47\) and \((\ell,k)\in\lbrace (167,11),(192,40)\rbrace\),
\item
\(M\simeq P_7-\#2;50\) and \((\ell,k)\in\lbrace (158,31),(171,41)\rbrace\),
\item
\(M\simeq P_7-\#2;52\) and \((\ell,k)\in\lbrace (161,32),(163,41)\rbrace\).
\end{enumerate}
The groups \(S\) are of order \(3^{20}\), class \(9\), coclass \(11\),
and their metabelianizations \(M\) are of order \(3^9\), class \(5\), coclass \(4\).
\end{lemma}


The lemmas have been proven with the aid of the
computational algebra system MAGMA
\cite{MAGMA},
developed at the University of Sydney
under supervision of John Cannon
\cite{BCP,BCFS}.

For the proof of Theorem
\ref{thm:ThreeStageTower33}
and Corollary
\ref{cor:ThreeStageTower33},
the invariants \(\varkappa(K)\), \(\tau^{(1)}(K)\) and \(\tau^{(2)}(K)\)
of \(3\)-class groups associated with algebraic number fields \(K\)
are first translated to invariants \(\varkappa(G)\), \(\tau^{(1)}(G)\) and \(\tau^{(2)}(G)\)
of finite \(3\)-groups, according to the Artin reciprocity law of class field theory
\cite{Ar1,Ar2}.

Then the group theoretic invariants are used as a search pattern
for Schur \(\sigma\)-groups, which are required as Galois group
\(\mathrm{Gal}(\mathrm{F}_3^\infty(K)/K)\)
of imaginary quadratic fields \(K\)
\cite{Ag,KoVe}.

Since Schur \(\sigma\)-groups have extremal root paths,
the search can be restricted to \(\sigma\)-groups connected by edges with maximal step size
equal to the nuclear rank of the respective parent.

The nuclear rank of the fork \(P_7\) is \(4\),
and the search pattern \(\varkappa(G)\), \(\tau^{(1)}(G)\) and \(\tau^{(2)}(G)\) is hit
by \(16\) descendants \(P_7-\#4;\ell\) for \(\varkappa(G)\sim (2114)\),
the eight in Lemma
\ref{lem:TypeF12}
and additionally \(\ell\in\lbrace 157,160,170,175,187,190,194,195\rbrace\), and
by \(16\) descendants \(P_7-\#4;\ell\) for \(\varkappa(G)\sim (2141)\),
the eight in Lemma
\ref{lem:TypeF13}
and additionally \(\ell\in\lbrace 112,115,122,132,135,137,141,147\rbrace\).

However, for the additional values of \(\ell\),
one component of \(\tau^{(2)}(G)\) uniformly evolutes
from \(\lbrack (27,9);(9,9,9,3),(27,3,3,3)^3\rbrack\) to \(\lbrack (27,9);(9,9,9,3),(81,3,3,3)^3\rbrack\)
at order \(3^{17}\), that is, for all descendants \(P_7-\#4;\ell-\#2;k-\#4;j\) with \(k=k(\ell)\) and \(1\le j\le 27\).
For the \(\ell\)-values in Lemma
\ref{lem:TypeF12}
and Lemma
\ref{lem:TypeF13},
no evolution occurs, and the pruned tree of \(\sigma\)-descendants is finite
and terminates in the Schur \(\sigma\)-groups
\(P_7-\#4;\ell-\#2;k-\#4;j-\#1;i-\#2;1\)
of order \(3^{20}\) as leaves. QED.


\section{Acknowledgement}
\label{s:Acknowledgement}

\noindent
The author gratefully acknowledges that his research was supported by the
Austrian Science Fund (FWF): P 26008-N25,
and by the Research Executive Agency of the European Union (EUREA).

Arithmetic invariants in the examples
and counterexamples
have been determined with the aid of class field theoretic MAGMA routines by Fieker
\cite{BCFS}.

Descendant trees of the forks \(P_7,P_8\) in the lemmas
have been constructed by means of the \(p\)-group generation algorithm by Newman
\cite{Nm}
and O'Brien
\cite{Ob},
which is also explained in
\cite{HEO}
and implemented in group theoretic MAGMA routines.




\end{document}